\documentstyle[12pt]{article}
\def\Bbb R{{\rm \bf R}}
\def\proclaim#1{\vskip2mm{\bf #1}\em}
\def\endproclaim{\em \vskip2mm}
\def\tag#1{\eqno(#1)}
\def\gathered{\begin{array}{c}}
\def\endgathered{\end{array}}
\def\text{\mbox}

\begin{document}

\title {The enclosure method for the heat equation}
\author{Masaru IKEHATA\\
Department of Mathematics,
Graduate School of Engineering\\
Gunma University, Kiryu 376-8515, JAPAN\\
Mishio KAWASHITA\\
Department of Mathematics,
Graduate School of Sciences\\
Hiroshima University, Higashi Hiroshima 739-8526, JAPAN}
\date{final}
\maketitle
\begin{abstract}
This paper shows how
the enclosure method which was originally introduced for elliptic equations
can be applied to inverse initial boundary value problems for parabolic equations.
For the purpose a prototype of inverse initial boundary value
problems whose governing equation is the heat equation is considered.
An explicit method to extract an approximation of the value of the support function at a given direction
of unknown discontinuity
embedded in a heat conductive body from the temperature for a suitable heat flux on the lateral boundary
for a fixed observation time is given.

\noindent
AMS: 35R30, 80A23

\noindent KEY WORDS: inverse initial boundary value problem, heat equation,
parabolic equation,
thermal imaging, cavity, corrosion, enclosure method
\end{abstract}

\section{Introduction}

The aim of this paper is to show how the {\it enclosure method}
which was originally introduced for elliptic equations in \cite{IE} can be applied to inverse initial
boundary value problems for parabolic equations in multi-dimensions.
We present only a simple case just for the explanation of the idea,
however, the same idea will work also for more general cases.

Let $\Omega$ be a bounded domain of $\Bbb R^m$, $m=2,3$ with a smooth boundary.
Let $D$ be an open subset of $\Omega$ with a smooth boundary and satisfy that:
$\overline D\subset\Omega$; $\Omega\setminus\overline D$ is connected.
We denote the unit outward normal vectors to $\partial\Omega$ and $\partial D$
by the same symbol $\nu$.
Let $T$ be an arbitrary {\it fixed} positive number.

Given $f=f(x,t),\,(x,t)\in\partial\Omega\times\,]0,\,T[$
let $u=u(x,t)$ be the solution of the initial boundary
value problem for the heat equation:
$$
\begin{array}{c}
\displaystyle \partial_tu-\triangle u=0\,\,\text{in}\,(\Omega\setminus\overline D)\times\,]0,\,T[,\\
\\
\displaystyle \frac{\partial u}{\partial\nu}=0\,\,\text{on}\,\partial D\times]0,\,T[,\\
\\
\displaystyle
\frac{\partial u}{\partial\nu}=f\,\,\text{on}\,\partial \Omega\times\,]0,\,T[,\\
\\
\displaystyle
u(x,0)=0\,\,\text{in}\,\Omega\setminus\overline D.
\end{array}
\tag {1.1}
$$

This paper considers the following problem for an explanation of the idea.

$\quad$

{\bf\noindent Inverse Problem.}  Assume that $D$ is {\it unknown}.
Extract information about the location and shape
of $D$ from the temperature $u(x,t)$ and
heat flux $f(x,t)$ for $(x,t)\in\,\partial\Omega\,\times\,]0,\,T[$.

$\quad$

The set $D$ is a model of the union of unknown {\it cavities} or subregions where
the heat conductivity is very low compared with that of the surrounded region $\Omega\setminus\overline D$.
Thus the problem is a
mathematical formulation of a typical inverse problem in thermal
imaging.  The method to solve this inverse problem may have possibility
of application in nondestructive evaluation.

There are extensive studies on the uniqueness and stability issues
on this kind of problems. See \cite{Ve5} and references therein
for the issues. In this paper we are concerned with the
reconstruction issue.  Note that in \cite{I4} Ikehata studied this type of
problems in one-space dimensional case by using the enclosure
method. In \cite{IC} he raised several questions related to the
method in multi-dimensional cases. However, it is still not clear
that the method works also for two or three-space dimensional
cases since the study in \cite{I4} fully makes use of the
speciality of one-space dimension.

Our main result is the following formula.

\proclaim{\noindent Theorem 1.1.}
Given $\omega\in S^{m-1}$ let $f$ be the function of $(x,t)\in\partial\Omega\times]0,\,T[$ having a parameter $\tau>0$
defined by the equation
$$\displaystyle
f(x,t;\tau)=\frac{\partial v}{\partial\nu}(x)\varphi(t),
\tag {1.2}
$$
where $v(x)=e^{\sqrt{\tau}x\cdot\omega}$ and $\varphi\in L^2(0,\,T)$
satisfying the condition:
there exists $\mu\in\Bbb R$ such
that
$$\displaystyle
\liminf_{\tau\longrightarrow\infty}
\tau^{\mu}
\left\vert\int_0^Te^{-\tau t}\varphi(t)dt\right\vert>0.
\tag {1.3}
$$
Let $u_f=u_f(x,t)$ be the weak solution of (1.1) for $f=f(x,t;\tau)$
and $h_D(\omega)=\sup_{x\in D}x\cdot\omega$.
Then the formula
$$\displaystyle
\lim_{\tau\longrightarrow\infty}
\frac{1}{2\sqrt{\tau}}
\log\left\vert
\int_{\partial\Omega}
\int_0^Te^{-\tau t}\left(v(x)f(x,t;\tau)-u_f(x,t)\frac{\partial v}{\partial\nu}(x)\right)dtdS\right\vert
=h_D(\omega),
\tag {1.4}
$$
is valid.

\endproclaim

\noindent
Note that: if $\varphi(t)$ is smooth on $[0,\,\delta[$ with $0<\delta<T$ and $t=0$ is not
a zero point with infinite order of $\varphi(t)$, then (1.3) is
satisfied for an appropriate $\mu>0$.

The function $h_D(\,\cdot\,)$ is called the {\it support function} of $D$ which gives
the signed distances from the origin of coordinates to the support plane ($m=3$), line ($m=2$)
of $D$.
Since the convex hull of $D$ is given by the formula
$\cap_{\omega\in S^{m-1}}\{x\in\Bbb R^m\,\vert\,x\cdot\omega<h_D(\omega)\}$, using the formula (1.4),
we obtain the set $\{x\in\Bbb R^m\,\vert\,x\cdot\omega<h_D(\omega)\}$ which gives an upper bound of the convex hull of $D$
from the direction $\omega$.

In \cite{I4} Ikehata considered the corresponding problem in one-space dimension for the equation $\partial_tu=(\gamma u_x)_x$
with a smooth $\gamma$.  Theorem 1.1 can be considered as an extension
of (2) of Theorem 2.5 in \cite{I4} in the case $\gamma=1$ to two and three-space dimensional cases.
Note that therein the parameter $\tau^2$ plays a role of $\tau$ in Theorem 1.1.

A brief outline of this paper is as follows.  Theorem 1.1 is
proved in Subsection 2.2. The proof is based on an integral
identity which is described in Subsubsection 2.2.1. Using the
identity, we give an asymptotic representation formula of the
integral
$$\displaystyle
\int_{\partial\Omega}\int_0^Te^{-\tau t}\left(v(x)f(x,t;\tau)-u_f(x,t)\frac{\partial v}{\partial\nu}(x)\right)dtdS
$$
together with an estimate of the leading term which is an integral
involving two Neumann-to-Dirichlet maps for the operator
$\triangle-\tau$. The estimate of the remainder term is described
in Subsubsection 2.2.2 as Lemma 2.1 and that of the leading term
is found in the same subsubsection as Lemma 2.2. To establish Lemma
2.1 we require an energy estimate of $u_f(\,\cdot\,,T)$: this together with the meaning of the weak solution of (1.1) is
found in Subsection 2.1.  In Section 3 we show how Theorem 1.1 can be extended also to extract
another information about $D$, that is the distance from a given $p$ outside $\Omega$ to $D$.

\section{The enclosure method}

In this section firstly we specify what we mean by the weak solution of
(1.1) and describe a necessary estimate.
Secondly by giving the proof of Theorem 1.1 we show how the enclosure method can be applied to the
heat equation.

\subsection{Preliminaries about the direct problem.}

We follow \cite{DL}.
For two separable Hilbert spaces $V$ and $H$ with $V\subset H\subset V'$, and a positive number
$T$, the space $W(0,T;V,V')$ is defined by $W(0,T;V,V')
=\{u\,\vert\,u\in L^2(0,T;V), u'\in L^2(0,T;V')\}$.  
Note that $u'$ means the derivative in $t\in]0,\,T[$.

Given $f\in L^2(0,\,T;H^{-1/2}(\partial\Omega))$ 
we say that $u\in W(0,\,T;H^1(\Omega\setminus\overline D),
(H^1(\Omega\setminus\overline D))')$ satisfy
$$\begin{array}{c}
\displaystyle
\partial_tu-\triangle u=0\,\,\text{in}\,(\Omega\setminus\overline D)\times]0,\,T[,\\
\\
\displaystyle
\frac{\partial u}{\partial\nu}=0\,\,\text{on}\,\partial D\times]0,\,T[,\\
\\
\displaystyle
\frac{\partial u}{\partial\nu}=f\,\,\text{on}\,\partial\Omega\times]0,\,T[
\end{array}
\tag {2.1}
$$
in the weak sense if the $u$ satisfies
$$
\displaystyle
\frac{d}{dt}\int_{\Omega\setminus\overline D}u(x,t)\varphi(x)dx
+\int_{\Omega\setminus\overline D}\nabla u(x,t)\cdot\nabla\varphi(x)dx
=<f(t),\varphi\vert_{\partial\Omega}>\,\text{in}\,(0,\,T),
\tag {2.2}
$$
in the sense of distribution on $(0,\,T)$ for all $\varphi\in H^1(\Omega\setminus\overline D)$.
By (1.50) on page 477 in \cite{DL} we have
$$\displaystyle
\frac{d}{dt}\int_{\Omega\setminus\overline D}u(x,t)\varphi(x)dx
=<u'(t),\varphi>
$$
in the sense of distribution on $(0,\,T)$ for all $\varphi\in H^1(\Omega\setminus\overline D)$.
Thus (2.2) is valid also for $t\in]0,\,T[$ a.e..

Note that by Theorem 1 on p.473 in \cite{DL}
we see that every $u\in W(0,\,T;H^1(\Omega\setminus\overline D), (H^1(\Omega\setminus\overline D))')$
is almost everywhere equal to
a continuous function of $[0,\,T]$ in $L^2(\Omega\setminus\overline D)$.  Further, we have:
$$\displaystyle
W(0,\,T;H^1(\Omega\setminus\overline D), (H^1(\Omega\setminus\overline D))')\hookrightarrow
C^0([0,\,T];L^2(\Omega\setminus\overline D)),
\tag {2.3}
$$
the space $C^0([0,\,T];L^2(\Omega\setminus\overline D))$ being equipped with the norm
of uniform convergence. Thus one can consider $u(0)$ and $u(T)$ as
elements of $L^2(\Omega\setminus\overline D)$. Then by Theorems 1 and 2 on p.512 and
513 in \cite{DL} we see that given $u_0\in L^2(\Omega\setminus\overline D)$ there
exists a unique $u$ such that $u$ satisfies (2.1) in the weak
sense and satisfies the initial condition $u(0)=u_0$.

Let $u_0=0$. By Remark 2 on p.512 and Theorem 3 on p.520 in
\cite{DL} we have the continuity of $u$ on $f$: there exists a
$C_T>0$ independent of $f$ such that
$$\displaystyle
\Vert u\Vert_{L^2(0,T;H^1(\Omega\setminus\overline D))}
\le C_T\Vert f\Vert_{L^2(0,T;H^{-1/2}(\partial\Omega))}.
\tag {2.4}
$$
Moreover, from (2.2) and (2.4) we have
$$\displaystyle
\Vert u'\Vert_{L^2(0,T;H^1(\Omega\setminus\overline D)')}
\le C_T\Vert f\Vert_{L^2(0,T;H^{-1/2}(\partial\Omega))}.
$$
This together with (2.3) and (2.4) yields one of the important estimates
in the enclosure method:
$$\displaystyle
\Vert u(T)\Vert_{L^2(\Omega\setminus\overline D)}
\le C_T\Vert f\Vert_{L^2(0,T;H^{-1/2}(\partial\Omega))}.
\tag {2.5}
$$
In the following subsection we denote by $u_f$ the weak solution
of (2.1) with $u(0)=0$ and this is the meaning of the weak
solution of (1.1).

\subsection{Proof of Theorem 1.1}

Define
$$\displaystyle
w_f(x;\tau)=\int_0^T e^{-\tau t}u_f(x,t)dt,\,\,x\in\Omega\setminus\overline D
$$
and
$$\displaystyle
g_f(x;\tau)=\int_0^Te^{-\tau t}f(x,t)dt,\,\,x\in\partial\Omega,
$$
where $\tau>0$ is a parameter.  This type of transform has been used
in the study \cite{I4} for the corresponding problem in a one-space
dimensional case.

\subsubsection{The basic identity}

The function $w_f$ satisfies
$$\begin{array}{c}
\displaystyle
(\triangle-\tau)w_f=e^{-\tau T}u_f(x,T)\,\,\text{in}\,\Omega\setminus\overline D,\\
\\
\displaystyle
\frac{\partial w_f}{\partial\nu}=0\,\,\text{on}\,\partial D,\\
\\
\displaystyle
\frac{\partial w_f}{\partial\nu}=g_f\,\,\text{on}\,\partial\Omega.
\end{array}
$$

Let $v=v(x)$ satisfy $(\triangle-\tau)v=0$ in $\Omega$.
Integration by parts yields
$$\displaystyle
\int_{\partial\Omega}\left(g_fv-w_f\frac{\partial v}{\partial\nu}\right)dS
=-\int_{\partial D}w_f\frac{\partial v}{\partial\nu}dS
+e^{-\tau T}\int_{\Omega\setminus\overline D}u_f(x,T)v(x)dx.
\tag {2.6}
$$
Let $p_f$ be the unique solution of the boundary value problem:
$$\begin{array}{c}
\displaystyle
(\triangle-\tau)p=0\,\,\text{in}\,\Omega\setminus\overline D,\\
\\
\displaystyle
\frac{\partial p}{\partial\nu}=0\,\,\text{on}\,\partial D,\\
\\
\displaystyle
\frac{\partial p}{\partial\nu}=g_f\,\,\text{on}\,\partial\Omega.
\end{array}
$$
Set $\epsilon_f=w_f-p_f$.
Since we have
$$\displaystyle
\int_{\partial\Omega}\left(g_fv-p_f\frac{\partial v}{\partial\nu}\right)dS=-\int_{\partial D}p_f\frac{\partial v}{\partial\nu}dS,
$$
from (2.6) we obtain the basic identity:
$$\begin{array}{c}
\displaystyle
\int_{\partial\Omega}\left(g_fv-w_f\frac{\partial v}{\partial\nu}\right)dS\\
\\
\displaystyle
=\int_{\partial\Omega}\left(g_fv-p_f\frac{\partial v}{\partial\nu}\right)dS
-\int_{\partial D}\epsilon_f\frac{\partial v}{\partial\nu}dS
+e^{-\tau T}\int_{\Omega\setminus\overline D}u_f(x,T)v(x)dx.
\end{array}
\tag {2.7}
$$
Note that $\epsilon_f$ satisfies
$$\begin{array}{c}
\displaystyle
(\triangle-\tau)\epsilon_f=e^{-\tau T}u_f(x,T)\,\,\text{in}\,\Omega\setminus\overline D,\\
\\
\displaystyle
\frac{\partial\epsilon_f}{\partial\nu}=0\,\,\text{on}\,\partial D,\\
\\
\displaystyle
\frac{\partial\epsilon_f}{\partial\nu}=0\,\,\text{on}\,\partial\Omega.
\end{array}
\tag {2.8}
$$
Multiplying the first equation in (2.8) by $\epsilon_f$ and integrating
over $\Omega\setminus\overline D$, we have
$$\displaystyle
\int_{\Omega\setminus\overline D}\vert\nabla\epsilon_f\vert^2dx+\tau\int_{\Omega\setminus\overline D}\vert\epsilon_f\vert^2dx
=-e^{-\tau T}\int_{\Omega\setminus\overline D}u_f(x,T)\epsilon_f(x)dx.
$$
This right hand side has the bound
$e^{-\tau T}\Vert u_f(T)\Vert_{L^2(\Omega\setminus\overline D)}\Vert\epsilon_f\Vert_{L^2(\Omega\setminus\overline D)}$.
Since $\tau>0$, we have immediately
$$\displaystyle
\Vert\epsilon_f\Vert_{H^1(\Omega\setminus\overline D)}
\le (\tau^{-1}+\tau^{-2})^{1/2}e^{-\tau T}\Vert u_f(T)\Vert_{L^2(\Omega\setminus\overline D)}.
\tag {2.9}
$$

\subsubsection{Two lemmas}

Now we choose a special $f$ having the form:
$$\displaystyle
f(x,t;\tau)=\frac{\partial v}{\partial\nu}(x)\varphi(t),
\,\,(x,t)\in\partial\Omega\times\,]0\,T[,
\tag {2.10}
$$
where $\varphi\in L^2(0,\,T)$.
For this $f$ we have
$$\displaystyle
\Vert f\Vert_{L^2(0,\,T;H^{-1/2}(\partial\Omega))}
=\Vert\varphi\Vert_{L^2(0,\,T)}\left\Vert\frac{\partial v}{\partial\nu}\right\Vert_{H^{-1/2}(\partial\Omega)}.
\tag {2.11}
$$

Since
$$\displaystyle
g_f(x;\tau)=\frac{\partial v}{\partial\nu}(x)\int_0^Te^{-\tau t}\varphi(t)dt,
$$
we have the expression of the first term of the right hand side of (2.7):
$$\displaystyle
\int_{\partial\Omega}\left(g_fv-p_f\frac{\partial v}{\partial\nu}\right)dS
=\int_0^Te^{-\tau t}\varphi(t)dt
\int_{\partial\Omega}\frac{\partial v}{\partial\nu}(R_{\emptyset}(\tau)-R_D(\tau))\frac{\partial v}{\partial\nu}dS,
\tag {2.12}
$$
where $R_{\emptyset}(\tau)$ and $R_D(\tau)$ are the Neumann-to-Dirichlet maps on $\partial\Omega$
for the operator $\triangle-\tau$ in $\Omega$ and $\Omega\setminus\overline D$ with the homogeneous Neumann
boundary condition on $\partial D$,
respectively.

Hereafter we choose a special $v$ in (2.10) which is a solution of the equation $(\triangle-\tau)v=0$ in $\Bbb R^m$:
$$\displaystyle
v(x)=e^{\sqrt{\tau}x\cdot\omega},
$$
where $\omega\in S^{m-1}$ and $\tau>0$.  Thus, we choose $f(x,t;\tau)$ as described in (1.2).

The following is an easy consequence of the estimates (2.5), (2.9)
and the growth order of $v$ as $\tau\longrightarrow\infty$.

\proclaim{\noindent Lemma 2.1.}
We have, as $\tau\longrightarrow\infty$
$$\displaystyle
-\int_{\partial D}\epsilon_f\frac{\partial v}{\partial\nu}dS
+e^{-\tau T}\int_{\Omega\setminus\overline D}u_f(x,T)v(x)dx
=O(\tau^{\gamma}e^{-\tau T}e^{2\sqrt{\tau}h_{\Omega}(\omega)})
\tag {2.13}
$$
where $\gamma$ is a positive constant and $h_{\Omega}(\omega)=\sup_{x\in\Omega}x\cdot\omega$.

\endproclaim
{\it\noindent Proof.}
First we show that the first term in (2.13) has the bound $O(\tau^{3/2}e^{-\tau T}e^{2\sqrt{\tau}h_{\Omega}(\omega)})$
as $\tau\longrightarrow\infty$.
Since $v$ satisfies $(\triangle-\tau)v=0$ in $\Omega$, we have, for all $\eta\in H^1(\Omega)$
$$\displaystyle
\int_{\partial\Omega}\frac{\partial v}{\partial\nu}\eta dS
=\int_{\Omega}\nabla v\cdot\nabla\eta dx+\tau\int_{\Omega}v\eta dx.
$$
This gives
$$\displaystyle
\left\vert\int_{\partial\Omega}\frac{\partial v}{\partial\nu}\eta dS\right\vert
\le(\Vert\nabla v\Vert_{L^2(\Omega)}+\tau\Vert v\Vert_{L^2(\Omega)})\Vert\eta\Vert_{H^1(\Omega)}.
$$
Since $\eta$ is arbitrary and the trace operator $H^1(\Omega)\longrightarrow H^{1/2}(\partial\Omega)$
has a bounded right inverse, we have
$$\displaystyle
\left\Vert\frac{\partial v}{\partial\nu}\right\Vert_{H^{-1/2}(\partial\Omega)}
\le C(\Vert\nabla v\Vert_{L^2(\Omega)}+\tau\Vert v\Vert_{L^2(\Omega)}).
\tag {2.14}
$$
Similarly we have
$$\displaystyle
\left\Vert\frac{\partial v}{\partial\nu}\right\Vert_{H^{-1/2}(\partial D)}
\le C(\Vert\nabla v\Vert_{L^2(D)}+\tau\Vert v\Vert_{L^2(D)}).
\tag {2.15}
$$
Since the trace operator $H^1(\Omega\setminus\overline D)\longrightarrow H^{1/2}(\partial D)$ is bounded,
we have
$$\displaystyle
\left\vert\int_{\partial D}\epsilon_f\frac{\partial v}{\partial\nu}dS\right\vert
\le C\left\Vert\frac{\partial v}{\partial\nu}\right\Vert_{H^{-1/2}(\partial D)}\Vert\epsilon_f\Vert_{H^1(\Omega\setminus\overline D)}.
$$
Now from this, (2.9) and (2.15) we obtain
$$\displaystyle
\left\vert\int_{\partial D}\epsilon_f\frac{\partial v}{\partial\nu}dS\right\vert
\le C(\tau^{-1}+\tau^{-2})^{1/2}(\Vert\nabla v\Vert_{L^2(\Omega)}+\tau\Vert v\Vert_{L^2(\Omega)})e^{-\tau T}
\Vert u_f(T)\Vert_{L^2(\Omega\setminus\overline D)}.
$$
This together with the special form of $v$, (2.5), (2.11) and (2.14) gives
the desired estimate.

It is clear that the second term of (2.13) has the bound $O(\tau e^{-\tau T}e^{2\sqrt{\tau}h_{\Omega}(\omega)})$
as $\tau\longrightarrow\infty$.

\noindent
$\Box$

The integral in (1.4) is just the left-hand side of (2.7).  This fact, (2.12) and Lemma 2.1
imply that
$$\begin{array}{c}
\displaystyle
\int_{\partial\Omega}
\int_0^Te^{-\tau t}\left(v(x)f(x,t;\tau)-u_f(x,t)\frac{\partial v}{\partial\nu}(x)\right)dtdS\\
\\
\displaystyle
=\int_0^Te^{-\tau t}\varphi(t)dt
\int_{\partial\Omega}\frac{\partial v}{\partial\nu}(R_{\emptyset}(\tau)-R_D(\tau))\frac{\partial v}{\partial\nu}dS
+O(\tau^{\gamma}e^{-\tau T}e^{2\sqrt{\tau}h_{\Omega}(\omega)}).
\end{array}
\tag {2.16}
$$
Thus, the main contribution on the limit (1.4) is the term given by two Neumann-to-Dirichlet maps.
For this term we need the following estimates from the both sides:

\proclaim{\noindent Lemma 2.2.}
There exist $\mu_1, \mu_2\in\Bbb R$, $\tau_0>0$, $C_1>0$ and $C_2>0$
such that, for all $\tau>\tau_0$
$$\displaystyle
C_1\tau^{\mu_1}\le
e^{-2\sqrt{\tau}h_D(\omega)}
\left\vert\int_{\partial\Omega}\frac{\partial v}{\partial\nu}(R_{\emptyset}(\tau)-R_D(\tau))\frac{\partial v}{\partial\nu}dS
\right\vert\le C_2\tau^{\mu_2}.
$$
\endproclaim

\noindent
Note that in Lemma 2.2 the lower estimate is essential
and the strict
values of $\mu_1$ and $\mu_2$ are not important to obtain formula (1.4).

From Lemma 2.2, (2.16) and (1.3), we have formula (1.4) in Theorem
1.1. More precisely, the condition (1.3) means that there exist
positive constants $C$ and $\tau_0$ such that, for all
$\tau\ge\tau_0$
$$\displaystyle
\tau^{\mu}\left\vert\int_0^T e^{-\tau t}\varphi(t)dt\right\vert\ge C.
$$
This together with (2.16) and the trivial estimate
$$\displaystyle
\left\vert\int_0^T e^{-\tau t}\varphi(t)dt\right\vert\le T^{1/2}\Vert\varphi\Vert_{L^2(0,\,T)}
$$
yields
$$\begin{array}{c}
\displaystyle
C\left\vert\int_{\partial\Omega}\frac{\partial v}{\partial\nu}(R_{\emptyset}(\tau)-R_D(\tau))
\frac{\partial v}{\partial\nu}dS\right\vert
+O(\tau^{\gamma+\mu}e^{-\tau T}e^{2\sqrt{\tau}h_{\Omega}(\omega)})\\
\\
\displaystyle
\le
\tau^{\mu}\left\vert\int_{\partial\Omega}
\int_0^Te^{-\tau t}\left(v(x)f(x,t;\tau)-u_f(x,t)\frac{\partial v}{\partial\nu}(x)\right)dtdS\right\vert\\
\\
\displaystyle
\le
T^{1/2}\Vert\varphi\Vert_{L^2(0,\,T)}\tau^{\mu}
\left\vert\int_{\partial\Omega}\frac{\partial v}{\partial\nu}(R_{\emptyset}
(\tau)-R_D(\tau))\frac{\partial v}{\partial\nu}dS\right\vert
+O(\tau^{\gamma+\mu}e^{-\tau T}e^{2\sqrt{\tau}h_{\Omega}(\omega)}).
\end{array}
$$
It follows from this and Lemma 2.2 that
$$\begin{array}{c}
\displaystyle
CC_1\tau^{\mu_1-\mu}\left(1+O(\tau^{\gamma+\mu-\mu_1}e^{-\tau T}e^{2\sqrt{\tau}(h_{\Omega}(\omega)-h_D(\omega))})\right)
\\
\\
\displaystyle
\le
e^{-2\sqrt{\tau}h_D(\omega)}\left\vert\int_{\partial\Omega}
\int_0^Te^{-\tau t}
\left(v(x)f(x,t;\tau)-u_f(x,t)\frac{\partial v}{\partial\nu}(x)\right)dtdS\right\vert\\
\\
\displaystyle
\le
C_2T^{1/2}\Vert\varphi\Vert_{L^2(0,\,T)}\tau^{\mu_2}\left(1+
O(\tau^{\gamma-\mu_2}e^{-\tau T}e^{2\sqrt{\tau}(h_{\Omega}(\omega)-h_D(\omega))})\right).
\end{array}
$$
Now the formula (1.4) follows from this.

The proof of Lemma 2.2 given below follows
the argument in the enclosure method applied to the Helmholtz
equation $(\triangle +k^2)u=0$ \cite{IE, IE2}.  It became
much simpler compared with that for the Helmholtz equation since our operator
$\triangle-\tau$ is strictly negative for $\tau>0$.

{\it\noindent Proof of Lemma 2.2.}
Consider the solution $R(x)$ of the following boundary value problem:
$$\begin{array}{c}
\displaystyle
(\triangle-\tau)R=0\,\,\text{in}\,\Omega\setminus\overline D,\\
\\
\displaystyle
\frac{\partial R}{\partial\nu}=\frac{\partial v}{\partial\nu}\,\,\text{on}\,\partial D,\\
\\
\displaystyle
\frac{\partial R}{\partial\nu}=0\,\,\text{on}\,\partial\Omega.
\end{array}
$$
Multiplying the first equation above by $R$ and integration over $\Omega\setminus\overline D$,
one gets
$$\displaystyle
\int_{\Omega\setminus\overline D}\vert\nabla R\vert^2dx+\tau\int_{\Omega\setminus\overline D}\vert R\vert^2dx
=-\int_{\partial D}\frac{\partial v}{\partial\nu}RdS.
\tag {2.17}
$$
Since the trace operator $H^1(\Omega\setminus\overline D)\longrightarrow H^{1/2}(\partial D)$ is bounded,
this right hand side has the bound
$$\displaystyle
\left\vert\int_{\partial D}\frac{\partial v}{\partial\nu}RdS
\right\vert
\le
2C\left\Vert\frac{\partial v}{\partial\nu}\right\Vert_{H^{-1/2}(\partial D)}\Vert R\Vert_{H^1(\Omega\setminus\overline D)}.
\tag {2.18}
$$
Now given $\eta>0$ let $\tau$ satisfy $\tau>\eta$.  Choose a $\epsilon>0$ in such a way that
$\min\,(1,\eta)>C\epsilon^2$.  Since we have
$$
\displaystyle
2C\left\Vert\frac{\partial v}{\partial\nu}\right\Vert_{H^{-1/2}(\partial D)}\Vert R\Vert_{H^1(\Omega\setminus\overline D)}
\le C\left(\epsilon^{-2}\left\Vert\frac{\partial v}{\partial\nu}\right\Vert_{H^{-1/2}(\partial D)}^2
+\epsilon^2\Vert R\Vert_{H^1(\Omega\setminus\overline D)}^2\right),
$$
it follows from (2.17) and (2.18) that
$$\displaystyle
\Vert R\Vert_{H^1(\Omega\setminus\overline D)}
\le\left(\frac{C}{\min\,(1,\eta)-C\epsilon^2}\right)^{1/2}\epsilon^{-1}
\left\Vert\frac{\partial v}{\partial\nu}\right\Vert_{H^{-1/2}(\partial D)}.
$$
Thus (2.15) implies that
$$\displaystyle
\Vert R\Vert_{H^1(\Omega\setminus\overline D)}\le C'(\Vert\nabla v\Vert_{L^2(D)}+\tau\Vert v\Vert_{L^2(D)})
\tag {2.19}
$$
where $C'>0$ is independent of $\tau>\eta$.
Noting that $R_{\emptyset}(\tau)(\partial v/\partial\nu)=v$, $R_D(\tau)(\partial v/\partial\nu)=v-R$ on $\partial\Omega$
and the functions $R$ and $v$ are real-valued, we have the integral identity
$$\begin{array}{c}
\displaystyle
\int_{\partial\Omega}\frac{\partial v}{\partial\nu}(R_{\emptyset}(\tau)-R_D(\tau))\frac{\partial v}{\partial\nu}dS\\
\\
\displaystyle
=\int_{\Omega\setminus\overline D}\vert\nabla R\vert^2dx+\tau\int_{\Omega\setminus\overline D}
\vert R\vert^2 dx
+\int_D\vert\nabla v\vert^2dx+\tau\int_D\vert v\vert^2dx.
\end{array}
$$
A combination of this, (2.19) and the speciality of $v$ gives
$$\displaystyle
2\tau\int_D e^{2\sqrt{\tau}x\cdot\omega}dx
\le \int_{\partial\Omega}\frac{\partial v}{\partial\nu}(R_{\emptyset}(\tau)-R_D(\tau))\frac{\partial v}{\partial\nu}dS
\le C_2\tau^3\int_D e^{2\sqrt{\tau}x\cdot\omega}dx.
$$
Now the conclusion follows from the corresponding estimate for the integral of $e^{2\sqrt{\tau}x\cdot\omega}$
over $D$ (cf. Propositions 3.1 and 3.2 in \cite{IE}).

\noindent
$\Box$

\section{The use of another $v$}

It is possible to use different $v$ from the
exponential solution.  In \cite{ISH} it is shown that: given $p\in\Bbb R^{m}\setminus\overline\Omega$
one can construct a $v\in H^2(\Omega)$ depending on $\tau>0$
and satisfying the equation $(\triangle-\tau)v=0$ in $\overline\Omega$
such that
$$\displaystyle
v(x)=e^{-\sqrt{\tau}\vert x-p\vert}
\left\{\vert x-p\vert^{-(m-1)/2}+O\left(\frac{1}{\sqrt{\tau}}\right)\right\}.
$$
Note that the leading term never vanish on $\overline\Omega$; in the case when $m=3$, one can
drop the term $O(1/\sqrt{\tau})$ in the $v$ above.

Using this $v$, we obtain the following formula.

\proclaim{\noindent Theorem 3.1.}
Let $p\in\Bbb R^m\setminus\overline\Omega$ and
replace $v$ of $f$ in (1.2) with the $v$ above.  Let $u_f=u_f(x,t)$ be the weak solution of (1.1) for this $f=f(x,t;\tau)$.
Then assuming (1.3), one has the formula
$$\displaystyle
\lim_{\tau\longrightarrow\infty}
\frac{1}{2\sqrt{\tau}}
\log\left\vert\int_{\partial\Omega}
\int_0^Te^{-\tau t}\left(v(x)f(x,t;\tau)-u_f(x,t)\frac{\partial v}{\partial\nu}(x)\right)dtdS\right\vert
=-d_D(p),
$$
where $d_D(p)$ denotes the distance from $p$ to $D$,
$$\displaystyle
d_D(p)=\inf\{\vert y-p\vert\,\vert\,y\in D\}.
$$

\endproclaim

Proceeding once more by the method used to prove (1.4), one knows that the key of the proof of
Theorem 3.1 is the following
lower estimate of the integral of $e^{-2\sqrt{\tau}\vert x-p\vert}$ over $D$.

\proclaim{\noindent Proposition 3.2.}
There exists $\mu\in\Bbb R$ such that
$$\displaystyle
\liminf_{\tau\longrightarrow\infty}\tau^{\mu}e^{2\sqrt{\tau}d_D(p)}\int_D e^{-2\sqrt{\tau}\vert x-p\vert}dx>0.
\tag {3.1}
$$

\endproclaim

{\it\noindent Proof.}
Choose $x_0\in\partial D$ such that $d_D(p)=\vert x_0-p\vert$.
Since $\partial D$ is smooth, one can find an open ball $B$ with $x_0\in\partial B$ and radius $r_0$
such that $B\subset D$.  Since the integrand is nonnegative, it suffices to prove (3.1) in the case when $D=B$.
We observe that
$$\displaystyle
\vert\nabla e^{-2\sqrt{\tau}\vert x-p\vert}\vert=2\sqrt{\tau}e^{-2\sqrt{\tau}\vert x-p\vert}.
$$
Using the co-area formula (Theorem 2.7.1. on page 76 in \cite{Z}), we obtain
$$\begin{array}{c}
\displaystyle
\int_B e^{-2\sqrt{\tau}\vert x-p\vert}dx
=\frac{1}{2\sqrt{\tau}}\int_{0}^{\infty}H^{m-1}(\{x\in B\,\vert\, e^{-2\sqrt{\tau}\vert x-p\vert}=t\})dt\\
\\
\displaystyle
=\int_{d_D(p)}^{\infty}H^{m-1}(S(s))e^{-2\sqrt{\tau}s} ds,
\end{array}
\tag {3.2}
$$
where $H^{m-1}$ denotes the $m-1$ dimensional Housdorff measure and
$S(s)=\{x\in B
\vert\,\vert x-p\vert=s\}$.

It is well known that $H^{m-1}$ agrees with the usual definition
of $m-1$-dimensional area on an $m-1$-dimensional $C^1$
submanifold of $\Bbb R^m$. This is a corollary of section 3.2 in
\cite{F} and consult also p.16 in \cite{Z}. Since $S(s)$ for
$d_D(p)<s<d_D(p)+r_0$ is a smooth surface ($m=3$), curve ($m=2)$,
$H^{m-1}(S(s))$ coincides with the $m-1$-dimensional area $\vert
S(s)\vert$. So the problem is reduced to give an estimate for
$\vert S(s)\vert$ from below as $s\downarrow d_D(p)$.

First consider the three-dimensional case.
The intersection of the set $\{x\in\,\Bbb R^3\,\vert\,\vert x-p\vert=s\}$ with $\partial B$ becomes a circle
on a plane and is given by the boundary of an open disc $S'(s)$ on the plane.
The radius $\delta(s)$ of $S'(s)$ is given by the equation
$$\displaystyle
\delta(s)^2
=r_0^2-\left\{\frac{r_0^2+(r_0+d_D(p))^2-s^2}{2(r_0+d_D(p))}\right\}^2
$$
and a simple computation yields $\delta(s)^2\ge C(s-d_D(p))$ as $s\downarrow d_D(p)$ for a positive constant
$C$.
This yields $\vert S'(s)\vert\ge\pi C(s-d_D(p))$ as $s\downarrow d_D(p)$.  Since $S'(s)$ is the projection
of $S(s)$ onto the plane mentioned above, we have $\vert S(s)\vert\ge\vert S'(s)\vert$.  Thus we obtain
$\vert S(s)\vert\ge\pi C(s-d_D(p))$ as $s\downarrow d_D(p)$.  From this and (3.2) we obtain the estimate
(3.1) with $\mu=1$.

Note that in the two-dimensional case we obtain $\vert S(s)\vert\ge 2C^{1/2}(s-d_D(p))^{1/2}$ 
as $s\downarrow d_D(p)$ and (3.1) is valid for $\mu=3/4$.

\noindent
$\Box$

\section{Conclusion and Remarks}

The procedure of extracting the support function of $D$ is extremely simple and summarized as follows.

\noindent
(i)  Give the direction $\omega\in S^{m-1}$.  Fix a large $\tau>0$ and give the heat flux across
$\partial\Omega$ over the time interval $]0,\,T[$:
$$\displaystyle
f(x,t;\tau)=\frac{\partial v}{\partial\nu}(x)
\varphi(t),\,\,(x,t)\in\partial\Omega\times]0,\,T[,
$$
where $v(x)=e^{\sqrt{\tau}x\cdot\omega}$ and $\varphi(t)$ satisfies (1.3) for a $\mu$.

\noindent
(ii) Measure the temperature $u_f(x,t)$ on $\partial\Omega$ over the time interval $]0,\,T[$.

\noindent
(iii)  Compute the quantity
$$\displaystyle
\frac{1}{2\sqrt{\tau}}
\log\left\vert\int_{\partial\Omega}
\int_0^Te^{-\tau t}\left(v(x)f(x,t;\tau)-u_f(x,t)\frac{\partial v}{\partial\nu}(x)\right)dtdS\right\vert
$$
as an approximation of $h_D(\omega)$.

As a corollary of Theorem 1.1 we have a {\it constructive} proof of the uniqueness of
recovering the convex hull of $D$ from infinitely many sets of the
temperature and heat flux on $\partial\Omega$ over the time
interval $]0,\,T[$.  Thus
our result can be considered as an extension of the enclosure method of {\it infinitely} many measurement version
\cite{IE}.

The method will cover also more general cases without serious
difficulty: inclusion, parabolic equations with
variable coefficients, Robin boundary condition.  And it may
be possible to apply our method to the corresponding problems for
hyperbolic equations and systems, Stokes system, some governing
equations in thermoelasticity, Maxwell systems, etc.  Such cases
will be reported in detail in forthcoming papers.

Note that in \cite{IK} we have already applied the enclosure method of
a {\it single} measurement version
\cite{IS} to the problem in three-dimensions.
It means that therein the heat flux $f$ is {\it fixed} and independent of large parameter $\tau$.
The information extracted therein is different from the information obtained in this paper.
The analysis developed therein is based on the potential theory and quite
delicate.  See also \cite{ISH} for an application of the enclosure method to an
inverse source problem for the heat equation in multi-space dimensions.

$$\quad$$

\centerline{{\bf Acknowledgements}}

MI was partially supported by Grant-in-Aid for
Scientific Research (C)(No. 18540160) of Japan  Society for
the Promotion of Science.
MK was partially supported by Grant-in-Aid for
Scientific Research (C)(No. 19540183) of Japan  Society for
the Promotion of Science.

$$\quad$$

\vskip1cm
\noindent
e-mail address

ikehata@math.sci.gunma-u.ac.jp

kawasita@math.sci.hiroshima-u.ac.jp
\end{document}